\documentclass[english]{article}
\usepackage[T1]{fontenc}
\usepackage[latin9]{inputenc}
\usepackage{amsmath}
\usepackage{amsthm}
\usepackage{amssymb}

\makeatletter
\theoremstyle{plain}
    \ifx\thechapter\undefined
	    \newtheorem{thm}{\protect\theoremname}
	  \else
      \newtheorem{thm}{\protect\theoremname}[chapter]
    \fi
\theoremstyle{plain}
\newtheorem{lem}{\protect\lemmaname}
\theoremstyle{remark}

\usepackage{geometry}

\makeatother

\usepackage{babel}
\providecommand{\lemmaname}{Lemma}
\providecommand{\remarkname}{Remark}
\providecommand{\theoremname}{Theorem}

\begin{document}
\title{On the theorem of Rubinstein}
\author{Rado\v s Baki\'{c}}
\date{}
\maketitle
\begin{abstract}

	Let $f(z)=\sum_{k=0}^{k=n}{n\choose k}a_{k}z^{k}$ and $r(z)=\sum_{k=n-p+1}^{k=n}\epsilon_{k}{n\choose k}a_{k}z^{k}$, with $|\epsilon_{k}|\le 1$ and $p<n-2$. Rubinstein proved the following theorem: if all zeros of $f(z)$ are in the region $|z|>R$, then all zeros of $f(z)+r(z)$ are in the region $|z|>\frac{R}{p+1}$. We give a new proof of this theorem that is more direct then the original proof. We prove that above theorem is also true under condition $|\epsilon_{k}|\le \frac{n}{e^2(p-1)}$ and $1<p\le n+1.$

\end{abstract}

\textbf{Key words} zeros of polynomial
\\

\\

Let $h_{1}(z)=\sum_{k=0}^{k=n}{n\choose k}a_{k}z^{k}$, and $h_{2}(z)=\sum_{k=0}^{k=n}{n\choose k}b_{k}z^{k}$ be two complex polynomials of degree $n$. Suppose also that their zeros are in the regions $|z|>r_{1}$, and $|z|>r_{2}$  respectively ($r_{1} r_{2}\neq 0$). Then well-known theorem of Szeg\"o [2] implies that zeros of their composite polynomial $h(z)=\sum_{k=0}^{k=n}{n\choose k}a_{k}b_{k}z^{k}$ are in the region $|z|>r_{1}r_{2}$. Let us note this is also true if degree of $h_{2}(z)$ is $k<n$. In that case polynomials $g_{1}(z)=z^nh_{1}(\frac{1}{z})$ and $g_{2}(z)=z^nh_{2}(\frac{1}{z})$ are both of degree n and  have zeros in the region $|z|<\frac{1}{r_{1}}$, and $|z|<\frac{1}{r_{2}}$ respectively. Then, again by the theorem of Szeg\"o, zeros of the composite polynomial of $g_{1}(z)$ and $g_{2}(z)$ are in the region $|z|<\frac{1}{r_{1}r_{2}},$ implying that zeros of $h(z)$ are in the region $|z|>r_{1}r_{2}$, as required. We shall use that fact in the further text.

As we said in the abstract,  our main result is:
 \begin{thm}
 \label{th1} 
	Let $f(z)=\sum_{k=0}^{k=n}{n\choose k}a_{k}z^{k}$ be a complex polynomial of degree $n$ such that all zeros are in the region $|z|>R$. Let also $r(z)=\sum_{k=n-p+1}^{k=n}\epsilon_{k}{n\choose k}a_{k}z^{k}$ with
	$|\epsilon_{k}|\le 1$ and $p<n-2$. Then, all zeros of $f_{1}(z)=f(z)+r(z)$ are in the region $|z|>\frac{R}{p+1}.$
	
\end{thm}

We shall now give a more direct proof of it. Case $p=1$ follows from the next lemma. Let us note that this lemma can be considered as a generalization of Corollary 1 of   [4].
\begin{lem}
	Suppose that $f(z)=\sum_{k=0}^{k=n}{n\choose k}a_{k}z^{k}$ has all zeros in the region $|z|>R$.  Then $f_{1}(z)=f(z)+\epsilon a_{n}z^{n}$ has all zeros in the region $|z|>\frac{R}{\sqrt[n]{|\epsilon |} +1  }$.
\end{lem}

\textbf{Proof:} By the well-known Coincidence theorem we have that $f(a)=a_{n}(a-c)^n$, for some complex $c$ depending on $a$, with $|c|>R$. If $a$ is zero of $f_{1}(z)$, then from $$0=f_{1}(a)=f(a)+\epsilon a_{n}a^n=a_{n}(a-c)^n +\epsilon  a_{n}a^n  $$
 follows easily that $|a|>\frac{R}{\sqrt[n]{|\epsilon |} +1  }$, which proves the lemma.

Let us now assume that $p>1$.  We shall use the following inequality, due to Biernacki [3]:

\begin{align}\label{ref1} 1+{n\choose 1}(p+1)+\cdots+{n\choose p-1}(p+1)^{p-1}<p^n,   
\end{align} for $1<p<n-2$.

\textbf{Proof of the Theorem 1:} Due to the Composition theorem of Szeg\"o (and our preliminary remarks) we can assume that in fact $f(z)=(1+z)^n$, even if the degree of $f_{1}(z)$ is less than $n$. So, we have to prove that zeros of the polynomial $f_{1}(z)=(1+z)^n+\sum_{k=n-p+1}^{k=n}{n\choose k}\epsilon_{k}z^{k}, \, |\epsilon_{k}|\le 1$, are all in the region $|z|>\frac{1}{p+1}.$ Suppose that it is not true. Then exist $b$, such that $f_{1}(b)=0$, and $|b|\le \frac{1}{p+1}$. From $f_{1}(b)=0$ it follows that $$ (1+b)^n=-\sum_{k=n-p+1}^{k=n}{n\choose k}\epsilon_{k}b^{k}.$$ Since $|1+b|\ge 1-|b|\ge 1-\frac{1}{p+1}=\frac{p}{p+1},$ we   conclude that $$
\left(\frac{p}{1+p} \right) ^n\le |1+b|^n =\left|  \sum_{k=n-p+1}^{k=n}{n\choose k}\epsilon_{k}b^{k}  \right|  \le \sum_{k=n-p+1}^{k=n}{n\choose k}\left( \frac{1}{p+1}\right) ^{k}, 
$$
i.e. $p^n\le \sum_{k=0}^{k=p-1}{n\choose k} (  p+1  ) ^{k}.$ This is a contradiction with \eqref{ref1}, and so the theorem is proved.
\\

In the original statement of Theorem 1 in [1], case $n=p-2$ is also included, but with incorrect proof. Proof of that case was based on the following inequality
\begin{align*} {q+2\choose q-1}(q+2)^{q-3}<q^{q+2},\,\text{for}\,q\ge 2
\end{align*}
which is false. In order to verify it, set $n=q+2$. Than we can rewrite our inequality into the following form
\begin{align*}
\frac{{n\choose 3}   }{ n^3 }<\left(1-\frac{2}{n} \right)^n,\, n\ge 4. 
\end{align*}
Taking limits on both sides we obtain
\begin{align*}
\frac{1}{6}\le \frac{1}{e^2}
\end{align*}
which  is a contradiction.
\\

For our next considerations we need the following lemma.
\begin{lem}
	 The following inequalities hold:
	 \begin{enumerate}
	 	\item $\frac{n}{e^2(p-1)}{n\choose p-1}(2+p)^{p-1}<p^n$, for $1<p\le n,\, 3\le n,$
	 	\item $\left( 1+{n\choose 1}(p+1)+\cdots+{n\choose p-1}(p+1)^{p-1}
	 	\right) \frac{n}{e^2(p-1)}<p^n$, for $1<p\le n+1,\, 1\le n.$
	 \end{enumerate}
\end{lem}
\textbf{Proof:}
\begin{enumerate}
	\item  Let $f(n,p)=\frac{ n{n\choose p-1}(2+p)^{p-1} }{ e^2(p-1)p^n  }.$ Then we have
	$$\frac{f(n,p)  }{f(n+1,p)}=\frac{ np(n-p+2)  }{ (n+1)^2 }    $$
	Condition $\frac{f(n,p)  }{f(n+1,p)}>1$ is equivalent to $n(n-(p-1))(p-1)>n+1$ and this is true, because at least one bracket on the left-hand side is greater than 1. Hence, $f(n,p)$ is decreasing on $n$ and in order to prove our inequality it is enough to prove that $ f(n,n+1)<1$, i.e.
	$$(1+\frac{2}{n+1})^{n+1}<e^2(1+\frac{2}{n+1})   $$ which is obviously true, because $ (1+\frac{2}{n+1})^{n+1}<e^2$.
	
	\item Case $n=p=2$ can be verified directly, and if $p=n+1$ then our inequality is again equivalent to $(1+\frac{2}{n+1})^{n+1}<e^2(1+\frac{2}{n+1}) $, as we had above. That means we can assume that $3\le n$ and $1<p\le n+1$. 
	In [3] Biernacki proved that 
	$$ 1+{n\choose 1}(p+1)+\cdots+{n\choose p-1}(p+1)^{p-1}<{n\choose p-1}(p+2)^{p-1}$$ for $1<p$.
	
	So, in order to prove our inequality it is sufficient to prove that $$ 
	\frac{ n  }{ e^2(p-1)    }{n\choose p-1} (2+p)^{p-1}<p^n
	$$ and this is already proved in 1.
	
\end{enumerate}
Our approach enables us also to obtain another theorem which is of the similar type as Theorem 1.
\\

\begin{thm}
	Let $ f(z)=\sum_{k=0}^{k=n}{n\choose k}a_{k}z^{k}$ be a complex polynomial of degree $n$, such that all zeros are in the region $|z|>R$. Let also $r(z)= \sum_{k=n-p+1}^{k=n}\epsilon_{k}{n\choose k}a_{k}z^{k}$, with $|\epsilon_{k}|\le \frac{ n  }{ e^2(p-1)  }$ and $1<p\le n+1$. Then zeros of $f_1(z)=f(z)+r(z)$ are all in the region $|z|>\frac{R}{p+1}$.
\end{thm}
\textbf{Proof:} Suppose that theorem is not true. Then exist $b$, such that $f_1(b)=0$, and $|b|\le \frac{1}{p+1}$. Using exactly same derivation as in the Theorem 1 we obtain that 
$$ p^n\le \sum_{k=0}^{k=p-1}{n\choose k}(p+1)^k \frac{ n  }{ e^2(p-1)  }. 
$$
This is a contradiction with  Lemma 2(part 2), therefore the theorem is proved.
\\
Let us note that Theorem 2 together with Lemma 1, enables us to change all coefficients of the polynomial, which is not the case with Theorem 1.
\\

Rado\v s Baki\'c\\
Teacher Education Faculty, Kraljice Natalije 43, Belgrade, Serbia\\
email: bakicr@gmail.com


\begin{thebibliography}{1}
\bibitem{lit1}\textsc{Rubinstein,  Z.}, \emph{Some inequalities for polynomials and their zeros},
Proc. Am. MAth. Soc., 16(1), 1965, 72-75.

\bibitem{lit2}\textsc{Szeg\"o, G.}, \emph{Bemerkungen zu einem Satz von J.H. Grace fiber die Wurzeln algebraisher Gleichungen}, Math. Z. 13, 1922, 28-55.

\bibitem{lit3}\textsc{Biernacki, M.}, \emph{Sur les zeros des polynomes}, Ann. Univ. Mariae Curie-Sklodowska Sect. A9, 1955, 81-98.

\bibitem{lit4}\textsc{Rahman, Q. I.}, \emph{The influence of coefficients on the zeros of polynomials}, J. London Math Soc. 36, 1961, 57-64.
\end{thebibliography}
\end{document}